\documentclass[10pt]{article}
\usepackage{amsmath}
\usepackage{amssymb}
\usepackage{mathrsfs}
\usepackage{amsfonts}
\usepackage{mathrsfs,amscd,amssymb,amsthm,amsmath,bm,graphicx}
\usepackage{graphicx}
\usepackage{multirow}
\usepackage{cite}
\usepackage{color}
\usepackage{enumitem}
\usepackage{array}
\usepackage{setspace}
\usepackage{float}
\usepackage{easyReview}
\usepackage{caption}
\usepackage[colorlinks,
linkcolor=black,
anchorcolor=black,
citecolor=black
]{hyperref}
\newcommand{\upcite}[1]{\textsuperscript{\textsuperscript{\cite{#1}}}}
\makeatletter
\renewcommand{\@seccntformat}[1]{{\csname the#1\endcsname}{\normalsize.}\hspace{.5em}}
\makeatother

\def \[{\begin{equation}}
\def \]{\end{equation}}

\newtheorem{thm}{Theorem}[section]

\newtheorem{lem}[thm]{Lemma}

\newtheorem{cor}[thm]{Corollary}

\begin{document}
\captionsetup[figure]{labelfont={bf},name={Fig.},labelsep=space}
\setlength{\baselineskip}{13pt}
\begin{center}{\Large \bf Analyses of Some Structural Properties on a Class of  Hierarchical Scale-free Networks
}

\vspace{4mm}

{\large Jia-Bao Liu $^{1, *}$, Yan Bao $^{1, *}$, Wu-Ting Zheng $^{1}$}\vspace{2mm}

{\small $^{1}$School of Mathematics and Physics, Anhui Jianzhu University, Hefei 230601, China}
\vspace{2mm}
\end{center}

\footnotetext{E-mail address: liujiabaoad@163.com, baoyanemail@163.com, zwtzjr@ahjzu.edu.cn.\\
}

\footnotetext{* Corresponding author.}

{\noindent{\bf Abstract.}\ \ Hierarchical networks actually have many applications in the real world. Firstly, we  propose a new class of  hierarchical networks with scale-free and fractal structure, which are the networks with triangles compared to traditional hierarchical networks. Secondly, we study the precise results of some structural properties to derive small-world effect and scale-free feature. Thirdly, it is found that the constructed network is sparse through the average degree and density. Fourthly, it is also demonstrated the degree distributions of hub nodes and the bottom nodes are the power law and exponential, respectively. Finally, we prove that clustering coefficient with a definite value $z$ tends to stabilize at a lower bound as $t$ iterates to a certain number, and the average distance of $G_{t}^{z}$ has a increasing relationship along with the value of $lnN_{t}$.

\noindent{\bf Keywords}:  Hierarchical Networks; Scale-free; Fractal Structure; Structural Properties\vspace{2mm}

\section{Introduction}\label{sct1}
Complex networks,\upcite{net1} such as information networks,\upcite{net2} social networks,\upcite{net3} biological networks,\upcite{net4}  is actually  simplified representations of  complex systems which reduces the system to an abstract structure that retains only basic connecting patterns features. In the past two decades, there are a large variety of researches  on complex networks emerging,\upcite{net5,net6,net7} which have attracted great attention from scientists and scholars. Topological properties and network dynamics of complex networks, the powerful tools to analyze networks, have two significant structural characteristics: small-world\upcite{Com} effect and scare-free\upcite{Bara01} feature, which take the complex theory to a  more superior standard and have a better application value. In recent studies,, researchers started  to derive some structural properties  of the complex networks such as  degree distribution,\upcite{Albert} cluster coefficient,\upcite{Zhang P} average path length,\upcite{Fron}  degree correlations\upcite{Pastor} and community structure\upcite{Clauset} with the purpose of  learning more about the theory behind the real world networks.

Except for the small-world effect and scale-free feature,  many realistic networks also have fractal structures and self-similarity, which is constructed  by  numerous continuous iterations.\upcite{fractal11,fractal17} Fractal complex networks are networks that are generated according to specific network generation rules,  making the networks infinitely more complex. For instance, Ref. \cite{Sierpinski} did a research on the Sierpinski pyramid network and derived its  power law strength.  Ref. \cite{Vicsek} studied  Vicsek fractal networks and calculate the mean geodesic distance.

In fact, there are  also many literatures about hierarchical networks. As for hierarchical fractal networks, Ref. \cite{Bara01} proposed a deterministic scale-free network which  was constructed with  fractal and hierarchical structures, showing that the degree sequence follows a power law. Since then, many researchers have studied various properties of this deterministic hierarchical network, and some researchers have also analysed the deformation of this network. Ref. \cite{Bara03} proved  that their deterministic hierarchical network model preserved a strict scaling law, and which had similar clustering behavior to many real complex networks. Then the various exactly solvable properties of the hierarchical network  were derived in Ref. \cite{Kazumoto}. Besides, Ref. \cite{zhangz} studied its average path length and did some analyses of deterministic network characteristics.

On the other hand, it noted that Ref. \cite{DAOHUA} generated a special hierarchical network having triangles whose root node connected other hub nodes of the same replicas and discussed some structural parameters of it.  At the same time,  combining with the promoted hierarchical networks in Ref. \cite{MNIU},  both networks are generated in an iterative manner. They gave us the inspiration to consider the method of generating the networks in this paper in a  special different way.

Inspired by the above literatures,  in the next work we propose a  new class of special hierarchical  networks with scale-free and fractal structure, which are the networks with triangles compared to traditional hierarchical networks. Then, we study the small-world effect and scale-free feature of the fractal networks and calculate the precise results of some property parameters, such as average degree, density, clustering coefficient, degree distribution and  average path length. At last, we give some deterministic conclusions about the networks.

The structure of the remained work is as follows.  In Sec. \ref{sct2}, the construction of the fractal networks is studied. In Sec. \ref{sct3}, several structural properties of the networks are computed and analysed. And conclusions are given in Sec. \ref{sct4}.

\section{Construction of the  fractal networks}\label{sct2}
In this section, we will generate a new class of special hierarchical fractal networks in an iterative way.
First of all, we denote the network as $G_{t}^{z}$ that changes with the number of iteration $t(t \in N, t\geq0)$ and the number of blocks $z(z \in N^{*}, z>2)$. Specifically, they are constructed in the following processes.
\\ \textbf{Step 1}. For $t=0$, we set a triangle as our initial graph denoted by $G_{0}$, having three nodes and edges. The vertex of the triangle is called  the hub node, getting a vertex set with an element denoted by $H_{0}$ in this way, and the other two nodes are called bottom nodes, which form a set as $B_{0}$.
\\ \textbf{Step 2}. For $t=1$, the graph $G_{1}$ can be obtained from $G_{0}$ in two manners:
\\(1) Have $z-1$ duplicates of $G_{0}$, so there are $z$ constituent blocks containing triangles in the network $G_{1}^{z}$. Then the largest hub node is the original hub node of $G_{0}$, and the sub-hub nodes are the hub nodes of the $z-1$ copies, and the bottom nodes of the duplicated triangles are called the last layer bottom nodes.
\\(2) Connect the hub node of $G_{0}$ to the last layer bottom nodes of the $z-1$ replicas. The set of the hub nodes and bottom nodes are named as $H_{1}^{z}$ and $B_{1}^{z}$, respectively.
\\ \textbf{Step 3}. Henceforth, for any time step $t>1$, the infinite network $G_{t}^{z}$ can be built from $z-1$ duplications of $G_{t-1}^{z}$ by connecting the first hub node of $G_{t-1}^{z}$ to the last layer bottom nodes of the copied duplications iteratively. In $G_{t}^{z}$, the original hub nodes of the $G_{t-1}^{z}$ and of its $z-1$ replicas express the new hub nodes of $G_{t}^{z}$, and the last layer bottom nodes of all $z-1$ replicas form the new generation of bottom nodes of $G_{t}^{z}$. The all hub nodes constitute the set $H_{t}^{z}$, and the bottom nodes of each copy form the set $B_{t}^{z}$. Fig. \ref{fig.1} illustrates the first three steps of the construction process of $G_{t}^{z}$ for $z=3$.

In the view of the growth  construction of the hierarchical fractal network $G_{t}^{z}$, we let $N_{t}$ be the total number of nodes, and $E_{t}$ be the total number of edges in $G_{t}^{z}$. When $t\geq0$, it is easy to get the expression of $N_{t}$ as
\begin{align}
 N_{t}= 3\cdot z^{t}. \label{N}
\end{align}
And we get the iterative expression of $E_{t}$ as
\begin{align}
E_{t}&= z\cdot E_{t-1}+2\cdot (z-1)^{t}.\label{e1}
\end{align}
With the initial condition $E_{0}=3$, we can deduce $E_{t}$ from Eq. (\ref{e1}) as
\begin{align}
E_{t}=(2z+1)\cdot z^{t}-2(z-1)^{t+1}.\label{E}
\end{align}

\begin{figure}[H]\
\centering\includegraphics[width=15cm,height=6cm]{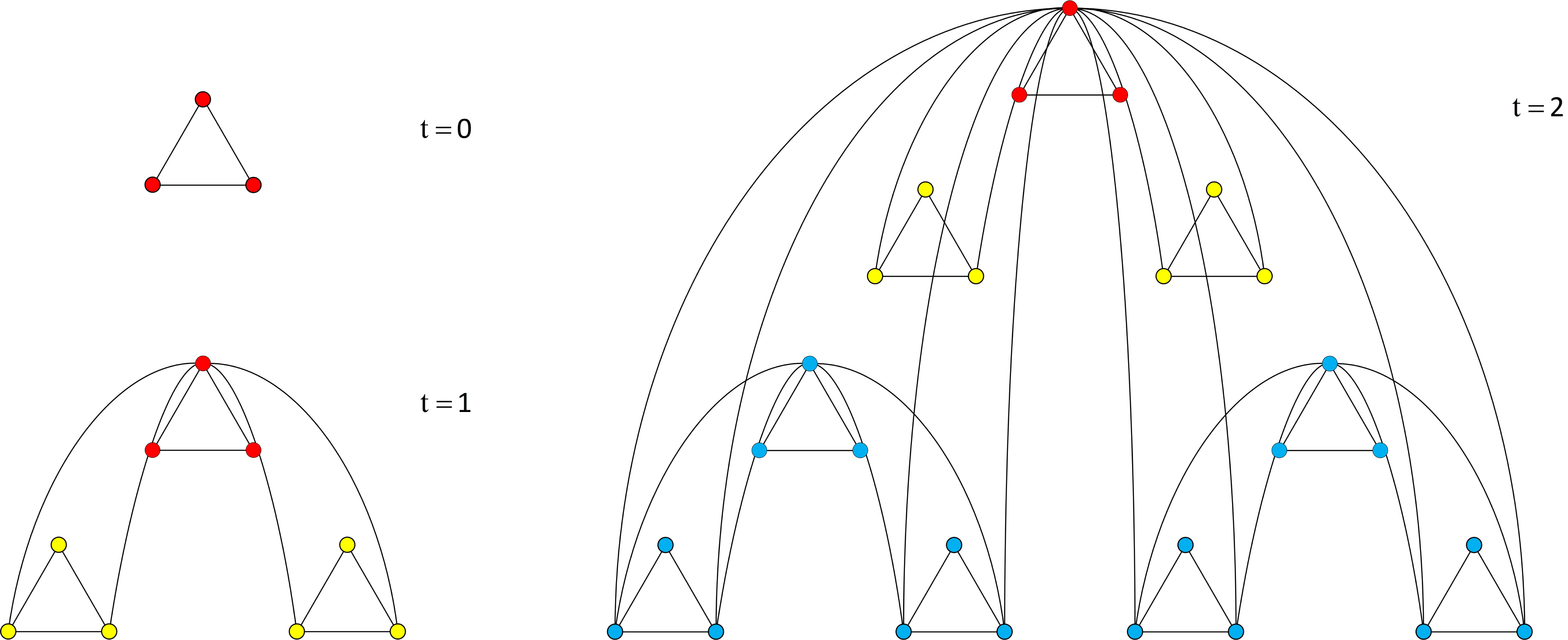}
\caption{The first two iterative processes of the networks for $z=3$.}
\label{fig.1}
\end{figure}

As we mentioned above,  the construction of the networks with hierarchical fractal structure have been finished.
In the next section, we will study some structural properties and discuss some features of the hierarchical fractal networks in detail.

\section{Analyses of some structural properties}\label{sct3}

In this section, we will calculate some important  structural property parameters of the $G_{t}^{z}$, these are, average degree, density, clustering coefficient, degree distribution and  average path length. At the same time, we will get the exact results of the structural  property parameters and derive the significance of different properties of the fractal networks. In addition, we discuss the some vital features of the fractal networks.

\subsection{Average degree} \label{sct3.1}
Average degree\upcite{MaF} is a simple structural parameter,  but it plays a significant role in complex network analysis, which is defined as the average degree of all nodes in the network, denoted by $\langle k_{t}\rangle$.  In general, it is  a indicator to determine whether a network is sparse.

Now, we aim to present the average degree of $G_{t}^{z}$ as follows.
\begin{thm} For\ any\ $t\in N, t\geq0\ and\ z\in N^{*}, z>2$,\ the\ average\ degree\ of\ the\ network\ $G_{t}^{z}$\  satisfies\
\begin{align*}
 \langle k_{t}\rangle & = \frac{4}{3}(1-z)(\frac{z-1}{z})^{t}+\frac{4z+2}{3}\\
 & \xrightarrow{t\rightarrow\infty} \frac{4z+2}{3}.
\end{align*}
\end{thm}

\noindent{\bf Proof.} According to Eq. (\ref{N}) and Eq. (\ref{E}), we can get $\langle k_{t}\rangle$ as Eq. (\ref{km}).
\begin{equation}\label{km}
\begin{aligned}
\langle k_{t}\rangle &=  \frac{2E_{t}}{N_{t}}  \\
&=\frac{1}{N_{t}} \sum_{i=1,j=1}^{N_{t}}(k_{i}+k_{j}) \\
&= -\frac{4}{3}z(\frac{z-1}{z})^{t+1}+\frac{4z+2}{3} \\ & \xrightarrow{t\rightarrow\infty} \frac{4z+2}{3},
\end{aligned}
\end{equation}
where $k_{i}$ and $k_{j}$ are degree of the hub nodes and bottom nodes respectively. \qed

In this way, when $t \rightarrow\infty $, the average degree tends to be a expression related to $z$. When $z$ is a certain value, the average degree is a constant. The average degree of most deterministic iterative networks is also a constant.\upcite{Bara01,Bara03} There are a small number of nodes in the network that have infinite connected edges, the number of nodes with few connected edges will increase infinitely. So it is clear to prove the network $G_{t}^{z}$ is a sparse network related to $z$.

\subsection{Density} \label{sct3.2}
Density\upcite{MaF} is also an important parameter to analyze network structure, which is defined as the ratio of the actual number of edges $M$ to the maximum possible number of edges containing $N$ nodes. In fact, it is also a indicator to determine whether a network is sparse.

Then, we can derive the density of $G_{t}^{z}$ as follows.
\begin{thm}For\ any\ $t\in N, t\geq0\ and\ z\in N^{*}, z>2$,\ the\ density\ of\ the\ network\ $G_{t}^{z}$\  has $$\rho \xrightarrow{t\rightarrow\infty}0.$$
\end{thm}

\noindent{\bf Proof.} For the network $G_{t}^{z}(t\in N, t\geq0, z\in N^{*}, z>2)$,  we can calculate the density of $G_{t}^{z}$ as Eq. (\ref{rho})
\begin{equation}\label{rho}
\begin{aligned}
\rho &=\frac{2E}{N(N-1)} \\
&= \frac{2\cdot((2z+1)\cdot z^{t}-2(z-1)^{z+1})}{3z^{t}\cdot(3\cdot z^{t}-1)} \\
& \xrightarrow{t\rightarrow\infty} 0.
\end{aligned}
\end{equation}
 When $t$ tends to $\infty$, the network density tends to $0$. This also further shows that the network $G_{t}$ is a sparse network. \qed
\subsection{Degree distribution} \label{sct3.3}
In complex network study, when a great number of connections appear on a small fraction of vertices, whereas the rest of the vertices have a small number of connections, the networks have a significant feature, that is, scale-free.\upcite{Bara19}  Many real world networks having scale-free feature  can be represented by degree distribution, which is defined as  the probability distribution of the degrees over the whole network in order to describe the scale-free feature of the network.

Next, we  will use statistical methods for depicting  degree distribution of the hub nodes and bottom nodes of $G_{t}^{z}$.

\begin{thm} $Let\ the\ degree\ of\  hub\ nodes\ and\ bottom\ nodes\  be\ denoted\ as\ k_{i}\ and\ k_{j}, \gamma = 1+ \frac{lnz}{ln(z-1)}, \theta= \frac{z-1}{z}.$ $The\ cumulative\ degree\ distribution\ of\ the\ network\ G_{t}^{z}\ is $
\begin{equation*}
P_{\mathrm{cum}}\left(K \geqslant k\right)=
\left\{\begin{array}{cc}
k_{i}^{1-\gamma}, & for\ hub\ nodes\ {i}, \\
\theta^{k_{j}}, & for\  bottom\ nodes\ {j}.
\end{array}\right.
\end{equation*}
\end{thm}

\noindent{\bf Proof.} In order to calculate the  degree distribution of $G_{t}^{z}$, we need to know the exact numbers   of  all nodes and degrees, which are also significant to get other structural properties. Therefore, let us firstly divide all the nodes and degrees into two categories: hub nodes set  $H_{t}^{z}$ and bottom nodes  set $B_{t}^{z}$, which is similar to the model in References \cite{Kazumoto, Bara01,MNIU}.

As for $H_{t}^{z}$, we can  know  that only one node  has the greatest degree $\frac{2(z-1)^{t+1}-2}{z-2}$ in  the $t$ iteration. And hub nodes of the $z$ replicas  all have the second greatest degree $\frac{2(z-1)^{t}-2}{z-2}$.  With such a classification, it is a fact that the $z$ newly generated replicas will not increase the degree of the hub nodes. So we can obtain that there are $(z-1)z^{t-i}$ nodes with degree $\frac{2(z-1)^{i}-2}{z-2}$.

As for $B_{t}^{z}$, we can find that there are $2(z-1)^{t}$ bottom nodes in the last layer having the greatest degree $t+2$ in  the $t$ iteration. In the next iteration, the number of degree of the original bottom nodes remains the same, and the number of degree of  bottom nodes in the two newly replicas  will add 1 due to the connecting with the largest hub node. Thus there are $2(z-1)^{j-2} z^{t-j+1}$ nodes with degree $j$.

Therefore, we get all the numbers of hub nodes $n(k_{i})$ and bottom nodes $n(k_{j})$ with their degree $k_{i}$ and $k_{j}$ in $G_{t}^{z}$ as shown in Table \ref{ta1}.

\begin{table}[H]\footnotesize
\setlength\tabcolsep{6pt}
\renewcommand{\arraystretch}{2}
\setlength{\abovecaptionskip}{0.05cm}
\centering \vspace{.3cm}
\caption{The numbers of hub nodes and bottom nodes with degree in $G_{t}^{z}$.} \label{ta1}
\begin{tabular}{m{1.8cm}<{\centering}m{1.5cm}<{\centering}m{1.5cm}<{\centering}m{1.5cm}<{\centering}m{1.6cm}
<{\centering}m{0.4cm}<{\centering}m{1.6cm}<{\centering}m{1.5cm}<{\centering}}
\hline
\multirow{2}{*}{hub nodes}
& $n(k_{i}) $    & $1$     & $z-1$      & $(z-1)\cdot z$         & $\cdots$ & $(z-1)\cdot z^{t-2}$  & $(z-1)\cdot z^{t-1}$
\\ \cline{2-8}
\multirow{2}{*}{} &$k_{i}$  &  $\frac{2(z-1)^{t+1}-2}{z-2}$     &  $\frac{2(z-1)^{t}-2}{z-2}$         & $\frac{2(z-1)^{t-1}-2}{z-2}$                                  & $\cdots$ & $2+2(z-1)$     & $2$
\\ \hline
\multirow{2}{*}{bottom nodes}
& $n(k_{j})$    & $2(z-1)^{t}$ & $2(z-1)^{t-1}$ & $2(z-1)^{t-2}z$ & $\cdots$ & $2(z-1)z^{t-2} $ & $2z^{t-1}$ \\ \cline{2-8}
\multirow{2}{*}{}  & $k_{j}$ & $t+2$   & $t+1$   & $t$      & $\cdots$ & $3$                                            & $2$ \\
\hline
\end{tabular}
\end{table}

In fact, Table \ref{ta1} is the degree sequence of $G_{t}^{z}$ in the jargon of graph theory, and it is the key condition for computing degree distribution. Combining the formula of cumulative degree distribution,\upcite{SN02} we can get the cumulative degree distribution of hub nodes and bottom nodes, respectively.

Taking the degree sequence of hub nodes into consideration, we have
\begin{equation}\label{pcmh}
\begin{aligned}
P_{cum}(k_{i})=&P(k_{i}\geq \frac{2(z-1)^{i}-2}{z-2})\\
=&\frac{z^{t-i+1}}{3z^{t}}\\
=&\frac{z}{3}z^{-i}\\
=&\frac{z}{3}((z-1)^{i})^{-\frac{lnz}{ln(z-1)}}. \\
\end{aligned}
\end{equation}
From Eq. (\ref{pcmh}), we can derive $P_{cum}(k)\varpropto k^{1-\gamma}$  with $\gamma = 1+\frac{lnz}{ln(z-1)}$ for hub nodes, so it is scale-free. When $z=3$, the vertices follow a power-law distribution with a power law of $1+\frac{ln3}{ln2}\approx 1.69$, and so on.

On the other hand, we obtain the degree distribution of bottom nodes in $G_{t}^{z}$ as
\begin{equation}\label{pcmb}
\begin{aligned}
P_{cum}(k_{j})=&P(k_{j}\geq j)\\
=&\frac{2(z-1)^{t}}{3z^{t}}\left(\frac{z}{z-1}\right)^{t-j+2} \\
=&\frac{2}{3}(\frac{z-1}{z})^{j-2}.\\
\end{aligned}
\end{equation}
According to Eq. (\ref{pcmb}), we can know $P_{cum}(k)\varpropto \theta^{k}$ with $\theta= \frac{z-1}{z}$ for bottom nodes. It indicates that the degree distribution of vertices is exponential  but not scale free. In a word, the hub nodes and bottom nodes have different scaling ways, so $G_{t}^{z}$ is a complex network which consists of hierarchical multi-fractal nature. \qed

\subsection{Clustering coefficient} \label{sct3.4}
As for another intriguing finding in complex network studies, clustering coefficient\upcite{Zhang P} is a property which can evaluate the small-world effect of the models. For example, the likelihood that two of your friends will be friends with each other in life reflects how connected your circle of friends is. The clustering coefficient can be used to quantify the probability that any two of your friends are friends with each other.

Theoretically,  Eq. (\ref{c1}) can abstractly describe the probability of connection among neighbors of any vertex $u$ with degree $k$.
\begin{align}
 C_{u}= \frac{n_{u}}{\frac{1}{2}k(k-1)}. \label{c1}
\end{align}
Then the  clustering coefficient denoted as $\overline{C}$ is defined as the average clustering coefficient of all nodes in $G_{t}^{z}$, just as Eq. (\ref{cm}).
\begin{align}
\overline{C}= \frac{\sum_{u\in G_{t}^{z}}C_{u}}{N_{t}}, \label{cm}
\end{align}
where $n_{u}$ expresses the existing edges between neighbors of any vertex $u$, $N_{t}$ is the total number of nodes in $G_{t}^{z}$. Then, we shall find a positive lower bound to verify the small-world property of $G_{t}^{z}$.
\begin{thm} $The\ clustering\ coefficient\ \overline{C_{t}}\  has\ a\  positive\ bound\ for\ different\ z\ blocks.\  When\ z=3$, $$\mathop{\lim}\limits_{t\rightarrow\infty} \overline{C_{t}} = 0.5279.$$
\end{thm}

\noindent{\bf Proof.} By definition, we can know that clustering coefficients $\overline{C_{t}}$ of the networks are not zero as there are triangles included in the network iteration process, which is different from the network in Ref. \cite{Kazumoto}.

For $t=0$, we find that the degree of every nodes, including one hub node and two bottom nodes, is 2. The every clustering coefficient is 1, so we get $\overline{C_{0}}=1$.

For $t=1$, the largest one hub node has coefficient with $\frac{2}{2z(2z-1)}$, and coefficient for $z-1$ sub-hub node is 1. The clustering coefficients of the $2(z-1)$ bottom nodes of the last layer, which all have a connecting edge to the largest hub node, is $\frac{2}{3\times2}$, and the other two bottom nodes have coefficients of 1. So we have
\begin{equation*}
\begin{aligned}
\overline{C_{1}}=&\frac{1}{3z}\left(\frac{1}{z(2z-1)}+\frac{2(z-1)}{3}+z+1\right)\\
=&\frac{5z+4}{9z^{2}(2z-1)}.
\end{aligned}
\end{equation*}

For any $t$, we firstly let $2(z-1)^{t+1}$ be $A^{t+1}$, so the original one hub node has the clustering coefficient of $\frac{2(z-2)^{2}}{(A^{t+1}-2)\times(A^{t+1}-z)}$. In the same way, $2(z-1)^{t}$ is $A^{t}$, and so on we can get $z-1$ hub nodes whose clustering coefficient are $\frac{2(z-2)^{2}}{(A^{t}-2)\times(A^{t}-z)}$,...; and there are $(z-1)z^{t-1}$ nodes with  clustering coefficients of $1$. Similarly, for the bottom nodes, there are $2(z-1)^{t}$ nodes with clustering coefficients of $\frac{2}{(t+2)\times(t+1)}$;  $2(z-1)^{t-1}$ nodes with clustering coefficients of $\frac{2}{(t+1)\times t}$;  $2(z-1)^{t-2} z$ nodes with clustering coefficients of $\frac{2}{t\times(t-1)}$,...; $2(z-1)z^{t-2}$ nodes with clustering coefficients of $\frac{2}{3\times2}$; $2z^{t-1}$ nodes with coefficients of 1. So the clustering coefficient $\overline{C_{t}}$ can be obtained as
\begin{equation}\label{C}
\begin{aligned}
\overline{C_{t}}
=&\frac{1}{3z^{t}}\left[\frac{2(z-2)^{2}}{(A^{t+1}-2)\times(A^{t+1}-z)}+\frac{4(z-1)^{t}}{(t+2)(t+1)}\right.\\
 &\left.
+\sum_{i=0}^{t-1}\frac{2(z-1)(z-2)^{2}z^{i}}{(A^{t-i}-2)\times(A^{t-i}-z)}
+\sum_{j=2}^{t+1}\frac{4(z-1)^{j-2}z^{t-j+1}}{j(j-1)}\right].
\end{aligned}
\end{equation}

For the sake of preliminarily estimate the relationship between the clustering coefficient $\overline{C_{t}}$ and $z$ and $t$, we make a mesh graph as plotted in Fig. \ref{fig.2}.
\begin{figure}[H]\
\centering\includegraphics[width=8.5cm,height=7cm]{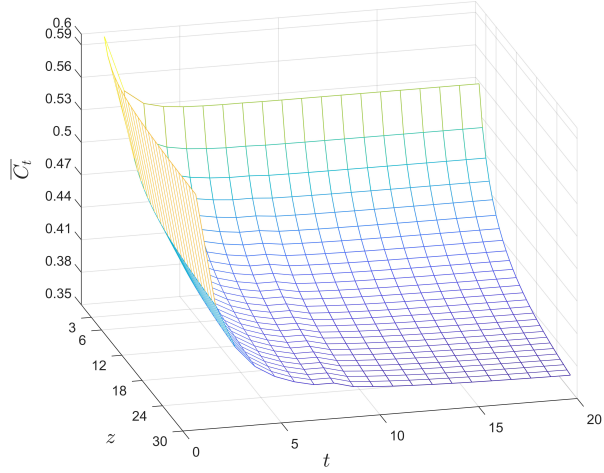}
\caption{The mesh graph about the relationship between $\overline{C_{t}}$ and $z$ and $t$.}
\label{fig.2}
\end{figure}
From Fig. \ref{fig.2}, we can consider that hierarchical networks with different $z$ have different clustering coefficients with the change of $t$ iteration times. We find that as $z$ increases, $\overline{C_{t}}$ tends to decrease, as does the law of $t$ and $\overline{C_{t}}$. In addition, when $z$ is a definite value, the value of $\overline{C_{t}}$ tends to stabilize with a lower bound as $t$ iterates to a certain number. Take $z=3$ as an example, as shown in Fig. \ref{fig.3} as follows.
\begin{figure}[H]\
\centering\includegraphics[width=8.5cm,height=7cm]{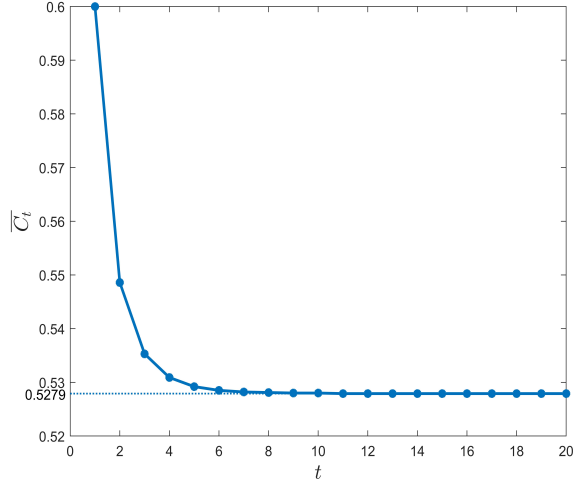}
\caption{The clustering coefficient of $G_{t}^{3}$. The dotted line shows the lower boundary of the curve.}
\label{fig.3}
\end{figure}
When $z=3$, according to Fig. \ref{fig.3} given above, we may conclude $\mathop{\lim}\limits_{t\rightarrow\infty} \overline{C_{t}} = 0.5279.$  In fact, at the time of the $9^{th}$ iteration, the change of the clustering coefficient values was already very little. So the clustering coefficient has a positive bound $0.5279$ in the network $G_{t}^{3}$. \qed

\subsection{Average path length} \label{sct3.5}
Average path length also called average distance,  is also known as a vital property of networks.\upcite{Fron} Then, we will derive the average distance of $G_{t}^{z}$ analytically. It can be calculated from the shortest distance $d_{uv}$ of nodes $u$ and $v$ for all possible pairs of nodes. When the iteration reaches the $t^{th}$ time, let $D_{t}$ be the total distance of all possible pairs of nodes in the network $G_{t}^{z}$ as Eq. (\ref{da1}).
\begin{equation}\label{da1}
D_{t}= \frac{1}{2}\sum_{u\neq v\in G_{t}^{z}} d_{uv}(t).
\end{equation}
Average path length denoted by $\overline{D_{t}}$ can be defined as Eq. (\ref{da2}).
\begin{equation}\label{da2}
\overline{D_{t}}= \frac{D_{t}}{N_{t}(N_{t}-1)/2}.
\end{equation}

The following work is necessary so as to calculate the average path of $G_{t}^{z}$. Firstly, we  consider the network $G_{t+1}^{z}$ as $z$ blocks abstractly,  which are named as $G_{t_{1}}^{z}$, $G_{t_{2}}^{z}, \cdots, G_{t_{z-1}}^{z},\ G_{t_{z}}^{z}$, respectively. We can see the abstract diagram of $G_{t+1}^{z}$  in Fig. \ref{fig.4}. For $G_{t+1}^{z}$, there are a total of $2(z-1)^{t+1}$ solid lines connecting the blocks, and each part has $2(z-1)^{t}$ connecting lines. Then we can know that the calculation of average path of  $G_{t+1}^{z}$ is to calculate the total distance  of pairs of nodes, which belong to the same self-similar structure $G_{t_{i}}^{z}(i=1, 2, \cdots, z)$ and belong to different branches, just as Eq. (\ref{D}).
\begin{figure}[H]\
\centering\includegraphics[width=13cm,height=7cm]{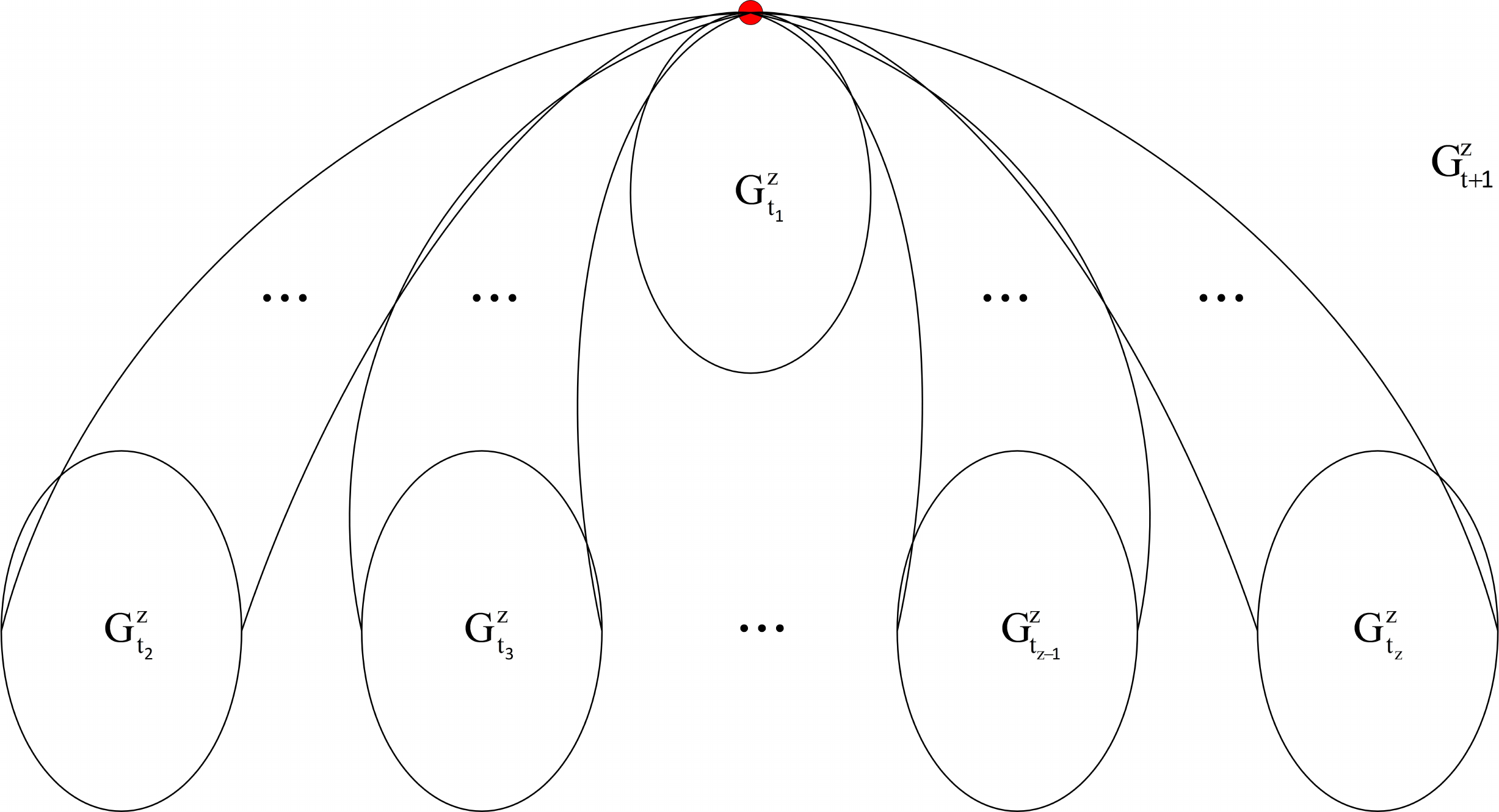}
\caption{The abstract diagram of $G_{t+1}^{z}$ with solid lines representing the lines between each component.}
\label{fig.4}
\end{figure}

\begin{equation}\label{D}
D_{t+1} = zD_{t}+\triangle_{t+1},
\end{equation}
where $\triangle_{t+1}$ represents the total distance of pairs of nodes belonging to different branches $G_{t_{i}}^{z}(i=1, 2, \cdots, z)$. Next,  we can divide the calculation process into two parts: (1) compute $\triangle_{t+1}$; (2) calculate  $D_{t}$, which are shown as following Lemma \ref{lem} and Theorem \ref{thd2}.
\begin{lem}\label{lem} When\ the\ network\ $G_{t}^{z}$\ with\ $z$\ blocks\ is\ in\ the\ $t^{th}$\ iteration,\
$$\triangle_{t+1}=6(z-1)z^{2t-1}\left(2z^{2}+(3t-1)z-3t\right).$$
\end{lem}

\noindent{\bf Proof.}
$\triangle_{t+1}$ can be expressed as some consisted components. Due to the self-similarity of the network, the specifics can be expressed as Eq. (\ref{trsum}).
\begin{equation}\label{trsum}
\begin{aligned}
\triangle_{t+1}
=& \triangle^{t_{1}^{z}, t_{2}^{z}}+ \triangle^{t_{1}^{z}, t_{3}^{z}}+\cdots+ \triangle^{t_{1}^{z-1}, t_{2}^{z}} \\
&+ \triangle^{t_{2}^{z}, t_{3}^{z}}+\triangle^{t_{2}^{z},t_{4}^{z}}+ \cdots+ \triangle^{t_{z-1}^{z}, t_{z}^{z}} \\
=& (z-1)\triangle^{t_{1}^{z}, t_{2}^{z}}+ C_{z-1}^{2}\triangle^{t_{2}^{z}, t_{3}^{z}},
\end{aligned}
\end{equation}
where $\triangle^{t_{1}^{z}, t_{2}^{z}}$ means that there are a pair of nodes $u$ and $v$ in $G_{t}^{z}$, $u$ is from $G_{t_{1}}^{z}$ and $v$ is from $G_{t_{2}}^{z}$, and other meanings can be followed by analogy.

Then, our task  is to calculate $\triangle^{t_{1}^{z}, t_{2}^{z}}$ and $\triangle^{t_{2}^{z}, t_{3}^{z}}$. Hence, we can have a arbitrary
node $u$  in $G_{t}^{z}$, and let $p_{u}(t)$ be the shortest distance from $u$ to the largest hub node, $q_{u}(t)$ be the shortest distance from $u$ to the last layer $2z^{t}$ bottom nodes. We also set $P_{t}$ stand the sum of $p_{u}(t)$ for all nodes in $G_{t}^{z}$ and make $Q_{t}$ express the total of $q_{u}(t)$. This way we can know the expression of $P_{t+1}$ as
\begin{align}
P_{t+1}=& \sum_{u \in G_{t_{1}}^{z}}p_{u}(t+1)+\sum_{u \in G_{t_{1}}^{z}}p_{u}(t+1)+ \cdots + \sum_{u \in G_{t_{z}}^{z}}p_{u}(t+1)  \nonumber \\
=& (z-1)\sum_{u \in G_{t}^{z}}\left(q_{u}(t)+1\right)+ \sum_{u \in G_{t}^{z}}p_{u}(t) \nonumber \\
=& (z-1)(Q_{t}+N_{t})+P_{t}. \label{eqp}
\end{align}
The expression of  $Q_{t+1}$ is
\begin{align}
Q_{t+1}=& \sum_{u \in G_{t_{1}}^{z}}q_{u}(t+1)+\sum_{u \in G_{t_{2}}^{z}}q_{u}(t+1)+\cdots + \sum_{u \in G_{t_{z}}^{z}}q_{u}(t+1) \nonumber\\
=& \sum_{u \in G_{t}^{z}}\left( p_{u}(t)+1\right)+ (z-1)\sum_{u \in G_{t}^{z}}q_{u}(t) \nonumber \\
=& P_{t}+N_{t}+(z-1)Q_{t}\label{eqq}.
\end{align}
Combining with $N_{t}=3z^{t}$, and the initial values of $P_{0}=2$ and $Q_{0}=1$, we can solve the Eq. (\ref{eqp})  for $P_{t}$ as
\begin{equation*}
P_{t}= z^{t-2}(4z^{2}+(6t-8)z-6t+6).
\end{equation*}
And  Eq. (\ref{eqq}) is solved for $Q_{t}$ as
\begin{equation*}
Q_{t}= z^{t-2}(z^{2}+(6t-2)z-6t+6).
\end{equation*}
With the solution of $P_{t}$ and $Q_{t}$, it is easy to  calculate $\triangle^{t_{1}^{z}, t_{2}^{z}}$ as
\begin{equation*}
\begin{aligned}
\triangle^{t_{1}^{z}, t_{2}^{z}}
=& \sum_{u\in G_{t_{1}}^{z},v\in G_{t_{2}}^{z}}d_{uv}(t+1) \\
=& \sum_{u\in G_{t_{1}}^{z},v\in G_{t_{2}}^{z}}\left( p_{u}+1+q_{v}\right)\\
=& \sum_{u\in G_{t_{1}}^{z}}\sum_{v\in G_{t_{2}}^{z}}p_{u}(t)+ \sum_{u\in G_{t_{1}}^{z}}\sum_{v\in G_{t_{2}}^{z}}\left(1+q_{v}(t)\right) \\
=& N_{t}\cdot P_{t} +N_{t}^{2}+N_{t}\cdot Q_{t} \\
=& 6z^{2t-2}(4z^{2}+(6t-5)z-6t+6).
\end{aligned}
\end{equation*}
And $\triangle^{t_{2}^{z}, t_{3}^{z}}$ is
\begin{equation*}
\begin{aligned}
\triangle^{t_{2}^{z}, t_{3}^{z}}
=& \sum_{u\in G_{t_{1}}^{z},v\in G_{t_{2}}^{z}}d_{uv}(t+1) \\
=& \sum_{u\in G_{t_{1}}^{z},v\in G_{t_{2}}^{z}}\left( q_{u}+1+1+q_{v}\right)\\
=& 2 \sum_{u\in G_{t_{1}}^{z}}\sum_{v\in G_{t_{2}}^{z}}\left(q_{u}(t)+1\right) \\
=& 2\left(N_{t}^{2}+N_{t}\cdot Q_{t}\right)\\
=& 12z^{2t-2}(2z^{2}+(3t-1)z-3t+3).
\end{aligned}
\end{equation*}
Therefore, Combining Eq. (\ref{trsum}), we can  get the specific solution of $\triangle_{t+1}$ as
\begin{equation}\label{tr}
\begin{aligned}
\triangle_{t+1}
=& (z-1)\triangle^{t_{1}^{z}, t_{2}^{z}}+ C_{z-1}^{2}\triangle^{t_{2}^{z}, t_{3}^{z}} \\
=& 6(z-1)z^{2t-1}\left(2z^{2}+(3t-1)z-3t\right).
\end{aligned}
\end{equation}
The proof is completed. \qed

Now, we can calculate the total distance $D_{t}$ in the following form.
\begin{thm}\label{thd2} When\ the\ network\ $G_{t}^{z}$\ with\ $z$\ blocks\ is\ in\ the\ $t^{th}$\ iteration,\
$$D_{t} = 3z^{t-2}\left(4z^{t+2}+(6t-8)z^{t+1}-6tz^{t}-3z^{2}+8z\right).$$
\end{thm}

\noindent{\bf Proof.} According to Eq. (\ref{D}) and Eq. (\ref{tr}), we can get a new iterative expression for $D_{t+1}$ as
\begin{equation*}
D_{t+1}= zD_{t}+6(z-1)z^{2t-1}\left(2z^{2}+(3t-1)z-3t\right).
\end{equation*}
Using the  initial condition $D_{0}=3$, we can get $D_{t}$ by iterative calculation  as
\begin{equation}\label{dd}
\begin{aligned}
D_{t}=& zD_{t-1}+ \triangle_{t} \\
=& z^{2}D_{t-2} +z \triangle_{t-1}+\triangle_{t} \\
=& z^{t}D_{0}+ z^{t-1}\triangle_{1}+ z^{t-2}\triangle_{2}+\cdots + \triangle_{t} \\
=& z^{t}D_{0} +\sum_{i=1}^{t}z^{t-i}\triangle_{i}.
\end{aligned}
\end{equation}
Then, we obtain
\begin{equation}\label{tri}
\sum_{i=1}^{t}z^{t-i}\triangle_{i} = 6z^{t-2}\left(2z^{t+2}+(3t-4)z^{t+1}-3tz^{t}-2z^{2}+4z\right).
\end{equation}
By taking Eq. (\ref{tri}) into Eq. (\ref{dd}), we get
\begin{equation}\label{dt}
D_{t} = 3z^{t-2}\left(4z^{t+2}+(6t-8)z^{t+1}-6tz^{t}-3z^{2}+8z\right).
\end{equation}
This completes the proof. \qed

Then, we can obtain theorem \ref{thD} of asymptotic formula as follows.
\begin{thm}\label{thD} For\ the\ network\ $G_{t}^{z},$\ one\ has\ $$\mathop{\lim}\limits_{t \rightarrow \infty}\frac{\overline{D_{t}}}{t}=\frac{4 (z-1)}{z^{2}}.$$
\end{thm}
\noindent{\bf Proof.}
In combination with $N_{t}=3z^{t}$ and Eq. (\ref{dt}), we can calculate
\begin{equation} \label{dm}
\begin{aligned}
\overline{D_{t}}=& \frac{D_{t}}{N_{t}(N_{t}-1)/2} \\
=& \frac{2\left(4z^{t+2}+(6t-8)z^{t+1}-6tz^{t}-3z^{2}+8z\right)}{z^{2}(3z^{t}-1)}.
\end{aligned}
\end{equation}
We can draw a three-dimensional mesh graph  of $\overline{D_{t}}$ as Fig. \ref{fig.5} to further explore the effects of $z$ and $t$ on $\overline{D_{t}}$.
\begin{figure}[H]\
\centering\includegraphics[width=8.5cm,height=7cm]{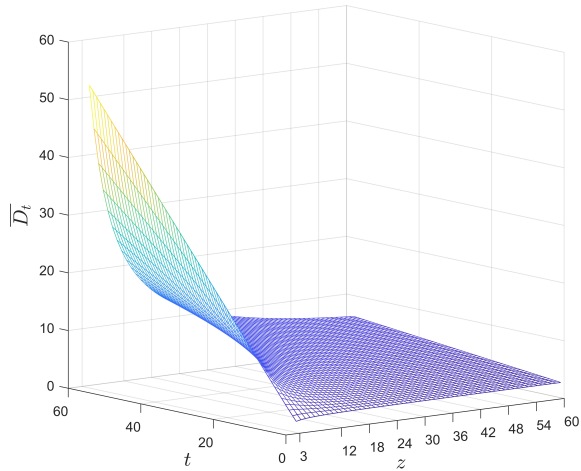}
\caption{The mesh graph of $\overline{D_{t}}$ for network $G_{t}^{z}$.}
\label{fig.5}
\end{figure}
It can be concluded from the mesh graph Fig. \ref{fig.5} that  $\overline{D_{t}}$ increases with the number of iterations $t$, and the increase of blocks $z$ will attenuate the growth rate of $\overline{D_{t}}$. Further, according to Eq. (\ref{dm}), we can figure out
\begin{equation}\label{lim}
\lim _{t \rightarrow \infty}\frac{\overline{D_{t}}}{t} =\frac{4 (z-1)}{z^{2}}.
\end{equation}

Therefore, it should be noted that the Eq. (\ref{lim}) of $\overline{D_{t}}$ has a linear increasing relationship with the number of iteration $t$ when $z$ has a certain value. We  can deduce $\overline{D_{t}}\varpropto ln N_{t}$ as $t$ has a logarithmic relationship with $N_{t}$. Moreover,  we find that the value of the asymptotic formula of $\overline{D_{t}}$ gradually decreases with the increase of $z$, and this result also verifies the phenomenon of Fig. \ref{fig.5}. \qed

Especially, when $z$ has a definite value for $z=3$, we will calculate the $\overline{D_{t}}$ of the network $G_{t}^{3}$ in Corollary \ref{corD}.
\begin{cor}\label{corD}
When $z=3$, we have average path length of the network $G_{t}^{3}$ as
$$\overline{D_{t}}=\frac{8lnN_{t}}{9ln3}+ O(3^{-t})\varpropto lnN_{t}.$$
\end{cor}

\noindent{\bf Proof.} From Eq. (\ref{N}), we  can get $N_{t}=3^{t+1}$ and $t=log_{3}N_{t}-1$. So we can continue to derive the relationship between $\overline{D_{t}}$  and $N_{t}$. Eq. (\ref{dm})  can be  rewritten as
\begin{equation*}
\begin{aligned}
\overline{D_{t}} =& \frac{8\cdot3^{t}(t+1)-2}{3^{t+2}-3} \\
=&\frac{8N_{t}log_{3}N_{t}-6}{9N_{t}-9}.
\end{aligned}
\end{equation*}
When the total number of nodes limit is  infinite $(N_{t}\rightarrow\infty)$, we find
\begin{equation}\label{D3}
\overline{D_{t}}=\frac{8lnN_{t}}{9ln3}+ O(3^{-t})\varpropto lnN_{t}.
\end{equation}
\begin{figure}[H]\
\centering\includegraphics[width=8.5cm,height=7cm]{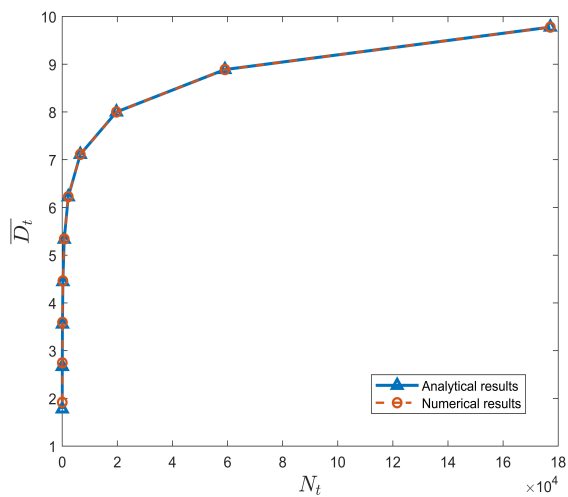}
\caption{The analytical and numerical results of  relationship between  $\overline{D_{t}}$ and $N_{t}$ for the network $G_{t}^{3}$.}
\label{fig.6}
\end{figure}
We perform rigorous calculation on the results of analytical simulation and find that the values are almost consistent with the theoretical  numerical results. Fig. \ref{fig.6} shows their comparison results. We find the average distance of $G_{t}^{3}$ increases along with the value of $lnN_{t}$ for the network. \qed
\\ \textbf{Remark \ref{sct3.5}.}
It is worth noting that the result of Theorem \ref{thD}, the asymptotic
formula, is consistent with  the unweighted hierarchical network in Ref. \cite{MNIU}. Also, our result of Corollary \ref{corD}, the analytical expression of $\overline{D_{t}}$ of $G_{t}^{3}$, is analogous with it researched by Ref. \cite{zhangz}.

Consequently, it  can be  proved that the existence of triangles in our network $G_{t}^{z}$ will not influence the asymptotic formula of $\overline{D_{t}}$, that is, have no  qualitative effect on the average path length. In addition, the average distance of $G_{t}^{z}$ has a increasing relationship with the value of $lnN_{t}$, where $N_{t}$ is the total number of nodes in $G_{t}^{z}$. Then, we can think the local structure of the network  $G_{t}^{z}$ has obvious character of collectivization, that is small-world effect.

\section{Conclusions}\label{sct4}
Through the above analyses, we study small-world effect and scale-free feature in the class of hierarchical networks with fractal structure.

In this paper, a new class of hierarchical  fractal networks are constructed iteratively. Under certain circumstances, the real network is also formed in the iterative way. Differing from the previously studied hierarchical networks, there are triangles in our researched  network $G_{t}^{z}$, whose the largest hub node has connecting edges to each bottom node of replicas in the last layer. Then, we obtain some results as follows.

(1) We study the average degree and density of the network so as to prove that the $G_{t}^{z}$ network is sparse.

(2) The degree distribution of hub nodes  and bottom nodes  in $G_{t}^{z}$ are demonstrated respectively,  where the hub nodes obey the power law distribution with a power index of $1+\frac{lnz}{ln(z-1)}$, and degree distribution of the bottom nodes is exponential.

(3) We explored the clustering coefficient of the network $G_{t}^{z}$ and determined  the  clustering coefficient with a definite value $z$ tends to stabilize at a lower bound as $t$ iterates to a certain number, such as  $\overline{C_{t}^{3}}$ is infinitely close to $0.5279$ as $t$ iterates to infinity.

(4) It is  found that the average distance of $G_{t}^{z}$  and the total number of nodes have a logarithmic growth relationship, that is, $\overline{D_{t}}\varpropto lnN_{t}$.

\section*{Funding}
Thanks for the support  of National Natural Science Foundation of China Grant 11601006.

\end{document}